\PassOptionsToPackage{dvipsnames}{xcolor}

\documentclass[10pt, journal, a4paper, twoside, twocolumn]{IEEEtran}
\usepackage[nolist]{acronym}
\usepackage[greek,english]{babel}

\usepackage{cite}

\usepackage[pdftex]{graphicx}

\graphicspath{{figures/}}
\DeclareGraphicsExtensions{.pdf,.png,.jpeg}

\usepackage[dvipsnames]{xcolor}
\definecolorset{cmyk}{}{}{%
    eth1, 1,    0.7,  0,   0.30;%
    eth2, 0.75, 0.4,  1,   0.40;%
    eth3, 1,    0.5,  0,   0;%
    eth4, 0.3,  0,    1,   0.55;%
    eth5, 0.2,  1,    0,   0.20;%
    eth6, 0,    0,    0,   0.7;%
    eth7, 0,    0.9,  0.8, 0.20;%
    eth8, 1,    0.25, 0.3, 0.1;%
    eth9, 0,    0.55, 1,   0.40;%
    eth10,0.6,  0,    1,   0 }

\definecolorset{RGB}{}{}{%
    eth1rgb,  31,  64, 122;%
    eth2rgb,  72,  90,  44;%
    eth3rgb,  18, 105, 176;%
    eth4rgb, 114, 121,  28;%
    eth5rgb, 145,   5, 106;%
    eth6rgb, 111, 111, 100;%
    eth7rgb, 168,  50,  45;%
    eth8rgb,   0, 122, 150;%
    eth9rgb, 149,  96,  19;%
    eth10rgb,140, 182,  60}

\usepackage[caption=false, font=footnotesize]{subfig}

\usepackage{amssymb}
\usepackage{amsmath}
\usepackage{mathrsfs} % for LaPlace / Fourier transform \mathscr{L}

\usepackage{url}
\usepackage[detect-all]{siunitx}
\DeclareSIUnit\pu{p.u.}
\DeclareSIUnit\voltampere{VA}
\DeclareSIUnit\var{var}

\usepackage{booktabs}
\usepackage{multirow}

\usepackage{bm} % robust bold face for math mode
\newcommand{\D}[2]{\mathrm{D}_{#2}#1}
\newcommand{\Dtwo}[2]{\mathrm{D}^2_{#2}#1}
\newcommand{\Dthree}[2]{\mathrm{D}^3_{#2}#1}

\newcommand{\Dn}[3][n]{\mathrm{D}^{#1}_{#3}#2}

\newcommand{\mat}[1]{\bm{#1}} % http://tex.stackexchange.com/questions/30619/what-is-the-best-symbol-for-vector-matrix-transpose
\renewcommand{\vec}[1]{\bm{#1}}

\newcommand{\revec}[1][i]{\vec{v}^{(#1)}}
\newcommand{\levec}[1][i]{\vec{u}^{(#1)}}
\newcommand{\levecT}[1][i]{\vec{u}^{(#1)\intercal}}
\newcommand{\eig}[1]{\lambda^{(#1)}}
\newcommand{\reeig}[1][i]{\sigma^{(#1)}}
\newcommand{\imeig}[1][i]{\omega^{(#1)}}
\newcommand{\damprat}[1][i]{\zeta^{(#1)}}

\newcommand{\res}[1][i]{\mat{r}^{(#1)}}
\newcommand{\kres}[1][i]{\frac{\res[#1]}{\eig{#1}}}

\newcommand{\trans}[1]{#1^\intercal}

\newcommand{\eqdot}{\quad .}
\newcommand{\eqcomma}{\quad ,}

\newcommand{\imag}{\mathrm{i}}
\newcommand{\eul}{\mathrm{e}}
\newcommand{\Rdim}[1]{\mathbb{R}^{#1}}

\newcommand{\vM}{\tilde{M}}
\newcommand{\vMi}{\tilde{M}_v}
\newcommand{\vK}{\tilde{K}}
\newcommand{\vKi}{\tilde{K}_v}
\newcommand{\vP}{{\overline{P}}}
\newcommand{\vPi}{{\overline{P}}_v}

\newcommand{\tp}[1][]{\mat{t}_\mathrm{S}^{#1}}

\newcommand{\Mp}[1][]{\mat{S}^{#1}}
\newcommand{\Mpnone}[1][]{\tilde{\mat{S}}^{#1}}
\newcommand{\Mpij}[1][]{(\mat{S})_{ij}^{#1}}
\newcommand{\Mpmax}{S_\infty}%{\lvert\Mp\rvert_\infty}
\newcommand{\Mpone}{S_1}%{\lvert\Mp\rvert_\infty}

\newcommand{\tROCOF}[1][]{\mat{t}_\mathrm{R}^{#1}}

\newcommand{\ROCOF}[1][]{\mat{R}^{#1}}
\newcommand{\ROCOFnone}[1][]{\tilde{\mat{R}}^{#1}}
\newcommand{\ROCOFij}[1][]{(\mat{R})_{ij}^{#1}}
\newcommand{\ROCOFmax}{R_\infty}%{\lvert\ROCOF\rvert_\infty}
\newcommand{\ROCOFone}{R_1}%{\normlone{\ROCOF}}

\newcommand{\czeta}{c^\zeta}
\newcommand{\zetamin}[1][]{\zeta^\mathrm{min}_{#1}}
\newcommand{\zetan}[1][]{\zeta_{\nu}^{(#1)}}
\newcommand{\zetanone}[1][]{\tilde{\zeta}_{\nu+1}^{(#1)}}

\newcommand{\cROCOF}{c^\mathrm{R}}
\newcommand{\cROCOFmean}{c^\mathrm{R1}}
\newcommand{\epsROCOF}[1][ij]{\epsilon^\mathrm{R}_{#1}}
\newcommand{\cepsROCOF}{c^{\epsilon\mathrm{R}}}

\newcommand{\cMp}{c^\mathrm{S}}

\newcommand{\Pbdg}{P^\mathrm{bdg}}
\newcommand{\Mbdg}{\tilde{M}^\mathrm{bdg}}

\newcommand{\hadam}{\odot}
\newcommand{\hadad}{\oslash}
\begin{document}

% TITEL
\title{\vspace*{18pt} On Placement of Synthetic Inertia\\ with Explicit Time-Domain Constraints}

\author{%
    \thanks{This research is supported by ETH Z{\"u}rich funds and by the European Union's Horizon 2020 research and innovation programme under grant agreement N$\circ$ 691800. This article reflects only the authors' views and the European Commission is not responsible for any use that may be made of the information it contains.}%
    \IEEEauthorblockN{Theodor S. Borsche and Florian D\"orfler\\}
    \IEEEauthorblockA{Automatic Control Lab, ETH Z\"urich\\
    borsche@alumni.ethz.ch, doerfler@control.ee.ethz.ch
  }%
}

%\author{Theodor Borsche, \IEEEmembership{Student Member, IEEE,} David J Hill, \IEEEmembership{Fellow, IEEE}%
%\thanks{Manuscript received May 35, 2064}%
%\thanks{All authors are with the Power Systems Laboratory at ETH Z\"urich, Z\"urich, Switzerland. Email: \{borsche\,\textbar\,ulbig\,\textbar\,andersson\}@eeh.ee.ethz.ch}%
%\thanks{Digital Object Identifier 12345}
%}

%\markboth{Journal of \LaTeX\ Class Files,~Vol.~6, No.~1, January~2007}%
%{Shell \MakeLowercase{\textit{et al.}}: Bare Demo of IEEEtran.cls for Journals}



\maketitle

\begin{acronym}[ENTSO-E]
    \acro{ROCOF}[RoCoF]{Rate of Change of Frequency}
    \acro{PSS}{Power System Stabilizer}
    \acro{PLL}{Phase-Locked Loop}
\end{acronym}

\begin{abstract}
Rotational inertia is stabilizing the frequency of electric power systems against small and large disturbances, but it is also the cause for oscillations between generators. As more and more conventional generators are replaced by renewable generation with little or no inertia, the dynamics of power systems will change. It has been proposed to add synthetic inertia to the power system to counteract these changes. This paper presents an algorithm to compute the optimal placement of synthetic inertia and damping in the system with respect to explicit time-domain constraints on the rate of change of frequency, the frequency overshoot after a step disturbance, and actuation input. A case study hints that the approach delivers reliable results, and it is scalable and applicable to realistic power system models.
\end{abstract}

%\begin{IEEEkeywords}
% IEEEtran, journal, \LaTeX, paper, template.
%\end{IEEEkeywords}
\acresetall
\section{Introduction}
Wind and solar generation have become some of the fasted growing energy sources world-wide --- initially for environmental reasons but increasingly also for economic aspects. This switch from conventional, thermal generation to renewables is also a switch from large synchronous generators to inverter-coupled generation. Photo-voltaic generation adds no inertia to power systems, and wind turbines in their most popular design add very little inertia. These developments have a serious effect on system dynamics \cite{Ulbig2014IFAC,tielens2016relevance}. Especially smaller interconnections are concerned about larger frequency incursions and \ac{ROCOF} after disturbances \cite{EirGrid2012, Ercot2013, statnett2016inertia}.

Inverters do not only decouple the inertia of wind turbines from the power system, but they can also be controlled to provide synthetic inertia and damping. This is achieved by adding a control loop that reacts to the \ac{ROCOF} or the more easily-measurable change in active power injection from the inverter \cite{Bevrani2014,DG-SB-BKP-FD:17}. Thus, synthetic inertia is becoming a design parameter of the power system, and some system operators even call for inertia-as-a-service\cite{EirGrid2012, Ercot2013}. 
We are hence faced with the questions of how much inertia we actually need, how to trade-off between virtual inertia and damping, where in the system has it the most beneficial effect, and how we can value the contribution of synthetic inertia? 

There are several approaches to answer these questions: Rakhshani et al. analyze the sensitivity of eigenmodes for tuning of virtual inertia and damping \cite{Rakhshani2016}. Poola et al. use $\mathcal H_2$ norms to minimize the energy content in the system frequency after a disturbance \cite{poolla2016placing,DG-SB-BKP-FD:17}. Pirani et al. use a related $\mathcal H_{\infty}$ criterion \cite{Pirani2017}.  Mesanovic et al. compare several approaches \cite{Mesanovic2016}. 
%\
While all of these approaches have their strengths and weaknesses, it is to be noted that they all optimize  objectives which are mere proxies for time-domain criteria such as \ac{ROCOF} or frequency deviations. Indeed, protective devices trigger based on the latter. For example, to avoid damage to generators due to vibrations and inadmissible  currents, the frequency deviation and \ac{ROCOF} must stay within limits. If these limits are violated, protection devices disconnect generators, likely starting a cascading failure. Moreover, none of the mentioned approaches can explicitly incorporate actuation constraints: the devices providing synthetic inertia and damping are limited in their power injection restricting their dynamic response.

 In this paper we extend our previous work  \cite{Borsche2015CDC} based on iterative eigenspace optimization with explicit time-domain constraints. Our approach considers not only a system-level objective specified in terms of eigenmodes but also explicit actuation constraints (power limits) as well as time-domain criteria on \ac{ROCOF} or frequency deviations.
 Our spectral performance criterion, the system damping ratio, and our optimization approach based on eigenspace sensitivities are similar to classic \ac{PSS} tuning for multi-machine systems\cite{Vournas1987}. Our approach is however more involved for two reasons. First, both the system dynamics as well as the input location are functions of the optimization parameters. %here not only the system dynamics matrix but also the input matrix is a function of the optimization parameter. 
 Second, we also optimize and enforce time-domain constraints which cannot be analytically found from the dynamic equations.
 Preliminary results \cite{BKP-DG-TB-SB-FD:17} suggest that the inertia distributions from the $\mathcal H_2$-based approach in \cite{poolla2016placing,DG-SB-BKP-FD:17} are very similar to the results obtained with the approach pursued in this paper,
 %. While \cite{poolla2016placing,DG-SB-BKP-FD:17} has a significantly lower computational effort, the scalability of our 
 but our approach seems to be more scalable which may be relevant for large systems. Additionally,  our approach can be used to optimize time-domain criteria while explicitly enforcing strict actuation constraints.

We illustrate our approach with a low-inertia version of the South-East Australian system adapted from Gibbard and Vowles \cite{Gibbard2014b}. This test-case has quite unique characteristics making it susceptible for low-inertia-driven instabilities: five areas are connected in a linear topology with  usually large flows from the outer regions to demand in the center of the system. Additionally, the Western and Northern ends are expected to see increase in wind and solar generation, respectively, thereby reducing the inertia in these loosely connected zones. Incidentally, a recent blackout was blamed on lacking fault-ride-through capabilities of wind farms in the western area \cite{AEMO2016}, albeit not on the lack of inertia.

The remainder of the paper is organized as follows: Section~\ref{sec:model} briefly introduces the modeling framework. Section~\ref{sec:vMvK_constraints} discusses to what extend a device can provide damping and inertia at the same time, which defines some constraints for our optimization. Section~\ref{sec:derivatives} analyzes the effect of synthetic inertia and damping on power system dynamics.  Section~\ref{sec:derivatives} describes how we compute the gradients of the non-linear placement problem. This is then used in Section~\ref{sec:optimization} to formulate an optimal inertia placement algorithm. Section~\ref{sec:testcase} gives a test case and showcases some results.

\section{Modeling} \label{sec:model}
For our design we use a small-signal model of a large power system. 
% Such models are commonly used for stability analysis with respect to damping of oscillatory modes and small disturbances, yet they are not ideal for large disturbance analysis as the linearization becomes increasingly inaccurate. However, the analysis in this paper would be intractable with a nonlinear system model, and while results should be checked against a nonlinear model we still expect to achieve meaningful insight using a linear model.
A power flow analysis of the system computes the steady state angles, voltages and active and reactive injections for a given load case. The dynamic model is then linearized around this steady state.

%\subsection{On notation}
%Derivatives are written in Euler's notation using the differential operator $\D{}{x}$, except for time derivatives in state space which use Newton's notation $\dot{x}$. Matrices and vectors are marked by bold letters, element-wise multiplication and division by $\hadam$ and $\hadad$. We follow conventions for eigenvalues $\eig{i}$, and left and right eigenvectors $\levec$, $\revec$.
%\begin{align}
 %   \levecT \mat{A} \revec &= \eig{i} = \reeig + \imag\,\imeig \\
  %  \levecT \revec &= 1 \eqdot \label{eq:evnormalization}
%\end{align}

\subsection{System dynamics}
The generator dynamics are modeled with six states describing the angle $\delta_i^\mathrm{G}$ and frequency $\omega_i^\mathrm{G}$ of the rotor, three fluxes in the machine and the excitation. In addition, each generator is equipped with an AVR, a PSS and a governor. Inputs to the system are the disturbances $\Delta P_k$,
%the generator mechanical power $P_{\mathrm{m}i}$ and voltage reference $V_{\mathrm{ref}i}$, 
and outputs are the rotor frequencies $\omega_i^\mathrm{G}$. 
Buses are connected via power lines, governed by the algebraic power flow equations.
%The linearized power flow equations are algebraic equations of the form
%\begin{equation}
%    P_i = -\sum\nolimits_{j \in \Omega_i} V_i^* V_j^* b_{ij} \cos(\delta_{ij}^*) (\delta_i - \delta_j)\eqcomma
%\end{equation}
%with $\Omega_i$ the set of neighboring buses, $V_i$ the bus voltage magnitude, $b_{ij}$ the inductance between buses, and the asterisk marking steady state values. 

The model so far only considers generator dynamics. Ignoring load dynamics may be inaccurate especially for grids with low inertia. It also means that there are no states related to load buses. % To place a controller or disturbance at a load bus, we need to recover these states.
In the spirit of \cite{Bergen1981}, we consider frequency-dependent load models and add motor loads to each bus to account for inertia in the system load. This recovers both the structure of the network interconnection and adds dynamic states %$\delta_i^\mathrm{L}$ and $\omega_i^\mathrm{L}$  for each load bus. 
for each load bus. 
%\fdmargin{Ich habe "$\delta_i^\mathrm{L}$ and $\omega_i^\mathrm{L}$" auskommentiert, da wir diese symbole nie wieder verwenden.} 
%With realistically chosen parameters, these motor loads do not significantly affect the system response to a disturbance. 
%
Kron reduction \cite{FD-FB:11d} is used to remove all remaining algebraic equations, and we arrive at the state-space model
\begin{subequations}
\label{sys:0}
\begin{align}
    \dot{\vec{x}}_0 &= \mat{A_0} \vec{x_0} + \mat{B_0} \vec{\Delta P} & \vec{\omega}^\mathrm{G} &= \mat{C_0} \vec{x_0}
\end{align}
\begin{align}
    \vec{x_0} & = 
        \trans{\begin{bmatrix} 
            \vec{\delta}^\mathrm{G} & \vec{\omega}^\mathrm{G} & \vec{\psi} & 
            \vec{x}^\mathrm{gen} & 
            \vec{\delta}^\mathrm{L} & \vec{\omega}^\mathrm{L}
        \end{bmatrix}}\\
%    \mat{A_0} &=
%        \begin{bmatrix}
%            \mat{0} & \mat{\omega_0} & \mat{0} & \mat{0} & \mat{0} & \mat{0} \\ 
%            -\inv{\mat{M}}\mat{L}^\mathrm{GG} & \mat{0} & \mat{0} & \mat{0} & \mat{0} & \mat{0} \\ 
%            & & & & & \mat{0} \\
%            & & & & & \mat{0} \\
%            \mat{0} & \mat{0} & \mat{0} & \mat{0} & \mat{0} & \mat{\omega_0} \\
%            \mat{0} & \mat{0} & \mat{0} & \mat{0} & -\inv{\mat{M}}\mat{L}^\mathrm{LL} & -\inv{\mat{M}}\mat{K}
%        \end{bmatrix} \\
    \mat{B_0} & = 
        \trans{\begin{bmatrix}
            \mat{0} & \mat{\frac{1}{M\,S_{\mathrm{B}i}}} & \mat{0} & \mat{0} & \mat{0} & \mat{0} 
        \end{bmatrix}} \\
    \mat{C_0} & = 
        \begin{bmatrix} 
            \mat{0} & \mat{\omega_0} & \mat{0} & \mat{0} & \mat{0} & \mat{0} 
        \end{bmatrix}
\end{align}
\end{subequations}
Our specific test case will be described in Section~\ref{sec:testsystem}.

\subsection{Synthetic inertia}
%\fdmargin{should we change notation from $\tilde K$ to $\tilde D$ as ``damping'' suggests ?}
Synthetic inertia can be provided by devices such as  batteries or supercaps \cite{tielens2016relevance}. On a system level we model synthetic inertia as a feedback loop of a grid-following converter. Each synthetic inertia block $v \in \mathcal{V}$ has bus frequency $\omega_v$ as input, and feeds power $\tilde{P}_v$ into the system according to the proportional-derivative (PD) control transfer function
\begin{align}
    \tilde{P}_v &= \frac{\vMi s + \vKi}{(T_{1v} s + 1)(T_{2v}s + 1)} \omega_v & v&\in \mathcal{V}\label{eq:virtualInertiaP}
\end{align}
We denote $\vMi$ as synthetic inertia, as it reacts proportional to the derivative of the frequency, and $\vKi$ synthetic damping, as it is proportional to frequency itself. The transfer function \eqref{eq:virtualInertiaP} has two poles -- one is needed for causality of the PD-control, the other accounts for the time constant of the \ac{PLL} to measure the frequency  $\omega_v$ at bus $v$.

A convenient state-pace representation of \eqref{eq:virtualInertiaP} is 
\begin{subequations}\label{eq:sysVi}%
\begin{align}%
    \mat{\tilde{A}_v} &= 
        \begin{bmatrix} 
            -\frac{T_{1v} + T_{2v}}{T_{1v} T_{2v}} & 1 \\[3pt]
            -\frac{1}{T_{1v} T_{2v}} & 0
        \end{bmatrix}
    & \mat{\tilde{B}_v} &= 
        \begin{bmatrix}
            \frac{\vMi}{T_{1v} T_{2v}} \\[3pt]
            \frac{\vKi}{T_{1v} T_{2v}}
        \end{bmatrix} \\
    %\mat{\tilde{C}_v} &= \begin{bmatrix} 1 & 0 \end{bmatrix} \label{eq:vM_feedback}
    \dot{\vec{\tilde{x}}}_v &=  \mat{\tilde{A}}_v \vec{\tilde{x}}_v + \mat{\tilde{B}}_v \omega_v 
    & \tilde{P}_v &= 
    %\underbrace{\begin{bmatrix} 1 & 0 \end{bmatrix}}_{\mat{\tilde{C}}_v} 
    \begin{bmatrix} 1 & 0 \end{bmatrix}\vec{\tilde{x}}_v
    =
    \mat{\tilde{C}}_v\vec{\tilde{x}}_v
\end{align}%
\end{subequations}%
The system \eqref{eq:sysVi} has two states $\tilde{x}$: the first state is power $\tilde{P}_v$ injected by the synthetic inertia device, and the second state is the measured frequency $\tilde{\omega}_v$. % (? unit is MW/s, not 1/s).
Connecting synthetic inertia \eqref{eq:sysVi} to the power system model \eqref{sys:0} gives the full system dynamics with disturbance $\vec{\Delta P}$ and output $\vec{\omega}^\mathrm{G}$ as
\begin{align}
    \begin{bmatrix} \dot{\vec{x}}_0 \\ \dot{\vec{\tilde{x}}} \end{bmatrix} &= 
    \underbrace{\begin{bmatrix}
        \mat{A_0} & -\mat{\Pi} \mat{\tilde{C}} \\ 
        \mat{\tilde{B}} \trans{\mat{\Pi}}  & \mat{\tilde{A}} 
    \end{bmatrix}}_{\displaystyle \mat{A}}
    \underbrace{\begin{bmatrix} \vec{x_0} \\ \vec{\tilde{x}} \end{bmatrix}}_{\displaystyle \vec{x}}+
    \underbrace{\begin{bmatrix} \mat{B_0} \\ \vec{0} \end{bmatrix}}_{\displaystyle \mat{B}} \vec{\Delta P} \\
    \vec{\omega}^\mathrm{G} &= \underbrace{\begin{bmatrix} \mat{C_0} & \vec{0} \end{bmatrix}}_{\displaystyle \mat{C}} \vec{x} \eqdot
\end{align}
The matrix $\mat{\Pi}$ with zero and unit entries maps the outputs $\tilde{P}_v$ and inputs $\omega_v$ of the synthetic inertia  \eqref{eq:sysVi} to system \eqref{sys:0}.  % = [\Pi_1, \ldots, \Pi_V]$; and $\mat{\Pi}_v^\$ mapping $\omega_v$ of \eqref{eq:virtualInertiaP} to $\mat{\tilde{B}}_v$, and $\mat{\Pi}_v$ mapping $\tilde{P}_v$ back to the input of $\mat{A}_0$. % The matrices $\mat{\tilde{A}}$, $\mat{\tilde{B}}$, $\mat{\tilde{C}}$ are block diagonal collections of \eqref{eq:sysVi}. The full model is of the form
% \begin{align}
%     \dot{x} &= \mat{A} {x} + \mat{B} \vec{u} & \vec{z} &= \mat{C} \vec{x}
% \end{align}
% with disturbances $\vec{u}$ and outputs $\vec{z}$, which are  $\Delta P$ and $\omega^\mathrm{G}$ respectively, confer \eqref{sys:0}.

%We will optimize over $\vM_i \in \mathcal{\tilde{M}}$ and $\vK_i \in \mathcal{\tilde{K}}$, which we collect in the set of parameters $\alpha \in \mathcal{\tilde{M}} \cup \mathcal{\tilde{K}}$. 
In the following, we will optimize over $\vMi$ and $\vKi$, which we collect in the set of parameters $\alpha \in \{\vMi, \vKi\}$. 
The derivative of $\mat{A}$ with respect to these parameters is
\begin{align}
    \D{\mat{A}}{\alpha} &= \begin{bmatrix} \mat{0} & \mat{0} \\ \D{\mat{\tilde{B}}}{\alpha} & \mat{0} \end{bmatrix} \\
    \D{\mat{\tilde{B}}}{\alpha} &= 
        \begin{cases}
            \trans{\begin{bmatrix} \frac{1}{T_{1v} T_{2v}} , 0 \end{bmatrix}} \trans{\mat{\Pi}_v} & \text{if~} \alpha \in \{\vMi\} \\[6pt]
            \trans{\begin{bmatrix} 0 , \frac{1}{T_{1v} T_{2v}} \end{bmatrix}} \trans{\mat{\Pi}_v} & \text{if~} \alpha = \{\vKi\} \\
        \end{cases} \eqdot
\end{align}
Observe that $\D{\mat{A}}{\alpha}$ is a sparse matrix with a single entry, making many of the following computations quite efficient.
\subsection{System response to disturbances}
The system response to disturbances can be described by two effects: 1) the damping ratio of oscillatory modes; and 2) the step response, e.g., after loss of generation.

%\textit{Damping ratio}
1) The \emph{damping ratio} $\damprat$ is obtained from the complex-conjugate eigenvalues $\eig{i} = \reeig + \imag\,\imeig$ of $\mat{A}$ as 
\begin{align}
  % \eig{i}_+ &= \reeig + \imag \imeig & \eig{i}_- &= \reeig - \imag \imeig \\
  \damprat &= \frac{-\reeig}{\sqrt{(\reeig)^2 + (\imeig)^2}} \eqcomma \label{eq:dampingratio}
\end{align}
and it is positive for stable and negative for unstable eigenvalues. Geometrically, $\damprat$ is the sine of the angle between the imaginary axis and a line from the origin to the eigenvalue. %$\zeta$ equals \num{+-1} for purely real eigenvalues.

%\textit{Step response} 
2) The \emph{step response} matrix $\mat{y}(\vec{t})$ can be computed at any $\vec{t}$ without explicit forward integration by 
\begin{align}
    \vec{y}(\vec{t}) &= \sum\nolimits\nolimits_i \kres \left( 1-\eul^{\eig{i}\vec{t}} \right)
    % y_{db}(t) &= \sum\nolimits\nolimits_i \kres \left( 1-\eul^{\eig{i}t} \right)
    \eqdot \label{eq:y_of_t}
\end{align}
This comes at the cost of computing the residues $\res$,
\begin{align}
    \res &= \mat{C} \revec \levecT \mat{B} \eqcomma \label{eq:residues}% & \kres &= \frac{\res}{\eig{i}} \eqcomma
\end{align}
which requires solving the eigenproblem of $\mat{A}$,
\begin{align}
    \levecT \mat{A} \revec &= \eig{i} = \reeig + \imag\,\imeig \\
    \levecT \revec &= 1 \eqcomma \label{eq:evnormalization}
\end{align}
where $\levec$, $\revec$ are the normalized left and right eigenvectors associated to the eigenvalue $\eig{i} = \reeig + \imag\,\imeig$.
The matrix $\mat{y}(\vec{t})$ collects the step responses at times $\vec{t}$ from any disturbance $\Delta P$ to any output $\omega^\mathrm{G}$. Note that $\vec{t}$ is a collection with the same dimension as $\vec{y}$, as we will later need the step response at different times for each disturbance-output pair.

%Equation \eqref{eq:y_of_t} yields the value of the step response $\mat{y}(t)$ at any $t$ without forward integration. This comes at the cost of solving the eigenproblem of $\mat{A}$ to find the residues $\res$. 
% {\color{eth4}If we investigate $D$ potential disturbances and observe $B$ states, the dimensions are $\mat{B} \in \Rdim{N\times D}$, $\mat{C} \in \Rdim{B \times N}$, and $\res$, $\mat{y}(t)$ are matrices $\in \Rdim{B \times D}$.}

We define the \emph{overshoot} $\Mp$ as the largest absolute value of the step response
%\begin{align}
%    \Mpi &:= \max\nolimits_{t} |y_{db}(t)| = |y_{db}(\tpi)| \eqcomma \label{eq:Mp}
%\end{align}
and the \emph{\ac{ROCOF}} $\ROCOF$ as the largest time derivative, occurring at $\tp$ and $\tROCOF$ respectively:
%\begin{align}
%    \ROCOFi &= \max\nolimits_{t} |\D{y_{db}}{t}(t)| = |\D{y_{db}}{t}(\tROCOFi)| \eqdot
%\end{align}
\begin{equation}
    \Mp = \mat{y}(\tp) = \max\nolimits_{\vec{t}} |\mat{y}|~, 
    \; \ROCOF = \mat{y}(\tROCOF) = \max\nolimits_{\vec{t}} |\D{\mat{y}}{t}|~. \label{eq:Mp} % (\vec{t})
\end{equation}
%Overshoot and \ac{ROCOF} occur at times $\tp$ and $\tROCOF$, for which $\mat{y}(\tp) = \Mp$ and $\mat{y}(\tROCOF)= \ROCOF$ holds. 
Finally, the largest overshoot and steepest \ac{ROCOF} are 
\begin{align}
    \Mpmax &:= \max\nolimits_{i,j}\, \Mp_{ij} \eqcomma & \ROCOFmax &:= \max\nolimits_{i,j}\, \ROCOF_{ij} \eqdot
\end{align}

\section{Constraints on synthetic inertia provision}\label{sec:vMvK_constraints}
Our placement algorithm considers a certain allowance of synthetic inertia $\vM$ and damping $\vK$, which in practice translates to a inverter-connected device, e.g., a battery or a supercap. On the timescales that we are considering, the power rather than the energy capacity of the inverter is limiting. The power injected by the device has to be lower than its power capacity $\bar{P}$, giving rise to the power constraint
\begin{align}
    |\vK \omega_k + \vM \dot{\omega}_k| & \leq \bar{P} &&\forall k \eqdot \label{eq:vKvM_leq_vPi}
\end{align}
In the following, we ask how much synthetic inertia and  damping one device can actually provide, and whether the provision of inertia restricts provision of damping? We answer these question by resorting to a data-driven approach.

\begin{figure}
    \centering
    \includegraphics{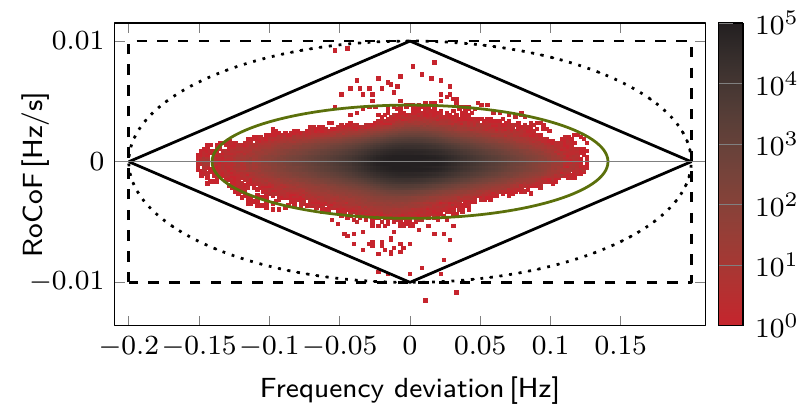}
    \caption{Scatter plot of frequency and \ac{ROCOF} data (one year of measurements  in \SI{1}{\second} resolution) in the CE interconnection (courtesy of RTE France). The color of each dot gives the frequency of each observation.}
    \label{fig:omegaOmegadot}
\end{figure}

%\subsection{Power requirements for inertia and damping provision}
%\textcolor{eth7}{[Habe den Abschnitt gek√ºrzt und lesbarer gemacht. Kann aber evtl ganz raus]}
%The power requested from an inertia device with $\vM_i$ and base power $\SBi$ for a given \ac{ROCOF} of $\dot{\omega}$ is
%\begin{equation}
%    \tilde{{P}}_i(\dot{\omega}) = -\frac{\vM_i \SBi}{\omega_0} \dot{{\omega}} \eqdot \label{eq:Pi}
%\end{equation}
%However, a change in frequency is the result of a given disturbance $\Delta P$. The power requested can than be approximated by
%\begin{equation}
%    \tilde{{P}}_i(\Delta P) = -\frac{\vM_i \SBi}{\omega_0} \frac{\omega_0}{M \SB} \Delta P = -\frac{\vM_i \SBi}{M \SB} \Delta P \eqdot \label{eq:Pi_system}
%\end{equation}
%When defining power capacity requirements, a system operator must either define a worst case \ac{ROCOF}, or a worst case fault and total inertia level.
%
%The same reasoning can be applied to damping provision, either assuming a worst case frequency deviation
%\begin{equation}
%    \tilde{P}_i = -\frac{\vK_i \SBi}{\omega_0} \omega \eqcomma
%\end{equation}
%or a worst case fault
%\begin{equation}
%    \tilde{P}_i = - \frac{\vK_i \SBi}{K \SB} \Delta P \eqdot
%\end{equation}

% \subsection{Trade-off between $\vM$ and $\vK$}
Figure~\ref{fig:omegaOmegadot} shows measurements from the continental European interconnection, covering one year in \SI{1}{\second} resolution. It is evident that the steepest \ac{ROCOF} and largest frequency deviation do not occur at the same time, and \SI{99.999}{\percent} of the data are contained in the green solid ellipse. 
The scaled 1-\nobreak, 2- and $\infty$-norm balls (black solid, dotted and dashed, respectively) are based on limits considered for normal system operation, namely \SI{200}{\milli\hertz} frequency deviation and \SI{10}{\milli\hertz\per\second} \ac{ROCOF}, and contain all but seven measurements.
% The black 1-norm is based on typical limits for frequency deviation and \ac{ROCOF}, and contains all but seven measurements. The more conservative 2-norm (dotted) and $\infty$-norm (dashed) also include almost all measurements.
%Intuitively, one does not expect to observe the worst case \ac{ROCOF} and worst case frequency deviation at the same time. 
The scaled 1- and 2-norm balls contain almost all observations, suggesting that simultaneous provision of $\vM$ and $\vK$ does not contradict itself. In the following we will show that the constraint for  $\vK$ and $\vM$ is the dual of the bounding norm in $(\omega,\dot{\omega})$ space.

We observe that the constraint \eqref{eq:vKvM_leq_vPi} is tight if $\omega_k$, $\dot{\omega}_k$ have the same sign. Henceforth, we will drop the absolute value. We  express $\vK$ as a function of $\vM$, $\omega$ and $\dot{\omega}$:
\begin{equation}
    \vK (\vM, \omega, \dot{\omega}) \leq \frac{\bar{P}}{\omega} - \frac{\vM \dot{\omega}}{\omega} \eqdot \label{eq:vKvM}
\end{equation}
Consider a set of observations $(\omega,\dot{\omega})$ contained in a scaled  $p$-norm ball of size $c>0$ and with scaling factor $h>0$:
\begin{align}
    \left( (h \omega_k)^p + \dot{\omega}^p_k\right)^{-p} & \leq c &&\forall k \eqdot  \label{eq:omegapnorm}
\end{align}
Assume the constraint  \eqref{eq:vKvM_leq_vPi} is tight for $(\omega,\dot{\omega})$ pairs on the boundary $(h\omega)^p + \dot{\omega}^p = c^p$. Hence, we can write \eqref{eq:vKvM} as
\begin{equation}
    \vK (\vM, \omega) \leq \frac{1}{\omega}\left(\vP - \vM \left(c^p - (h\omega)^p\right)^{-p}\right) \eqdot \label{eq:vKvM2}
\end{equation}
Thus, for each $\omega$ there is a linear upper bound \eqref{eq:vKvM2} in $(\vM,\vK)$ space, constraining the choice of $\vK$ depending on $\vM$. The limiting constraint on $\vK$ for a given $\vM$ can be found by
minimizing the right-hand side of \eqref{eq:vKvM2} with respect to $\omega$. Aside from critical points at $|\omega| = \infty$, we obtain the others by setting the derivative of the right-hand side of \eqref{eq:vKvM2} to zero
\begin{equation}
-\vP + (c^p - (h\omega)^p)^{q} c^p \vM = 0 \label{eq:minK}
\end{equation}
%\begin{multline}
%    \min\nolimits_\omega \vK(\vM, \omega)  \Leftrightarrow \D{\vK}{\omega}|_{\vM} = 0\\
%    \Leftrightarrow  -\vP + (c^p - \omega^p)^{q} c^p \vM = 0 \label{eq:minK}
%    % -\frac{\vPi}{\omega^2} + (c^p - \omega^p)^q\frac{c^p \vM} { \omega^2} = 0 \label{eq:minK}
%\end{multline}
with $q = \frac{p}{p-1}$. 
After solving \eqref{eq:minK} for $\omega$, substituting the solution in \eqref{eq:vKvM}, and after some reformulations, we arrive at
\begin{equation}
    \left(\frac{\vK}{h}\right)^q + \vM^q = \left(\frac{\vP}{c}\right)^q \eqdot \label{eq:qnorm}
\end{equation}
Hence, $\vK$ and $\vM$ are bounded by the dual norm of \eqref{eq:omegapnorm}.

%{\color{eth6}
%Equation \eqref{eq:minK} is fulfilled for 
%\begin{align}
%    \omega^p = c^p - \left(\frac{c^p \vM}{\vPi}\right)^q 
%    \quad \overset{\eqref{eq:omegapnorm}}{\Leftrightarrow} \quad 
%    \dot{\omega}^p =\left(\frac{c^p \vM}{\vPi}\right)^q \label{eq:omegaplimiting}
%\end{align}
%Inserting \eqref{eq:omegaplimiting} in \eqref{eq:vKvM} yields
%\begin{align}
%    \vK \omega &= \vK \left(c^p - \left(\frac{c^p \vM}{\vPi}\right)^q\right)^\frac{1}{p} \label{eq:vKomega1}\\
%    &\overset{\eqref{eq:omegaplimiting}}{=} \vPi - \vM \left(\frac{c^p \vM}{\vPi}\right)^\frac{1}{p-1} \\
%    &= \frac{\vPi}{c^p} \left(c^p - \left(\frac{c^p \vM}{\vPi}\right)^q\right)  \label{eq:vKomega2}
%\end{align}
%Combining \eqref{eq:vKomega1} and \eqref{eq:vKomega2} gives
%\begin{multline}
%     \vK = \frac{\vPi}{c^p} \left(c^p - \left(\frac{c^p \vM}{\vPi}\right)^q\right)^{-q} \\
%     \Leftrightarrow \vK^q + \vM^q = \left(\frac{\vPi}{c}\right)^q
%\end{multline}
%which shows that $\vK$ and $\vM$ are bounded by a ${q}$-norm, which is the dual of a $p$-norm since $q = \frac{p}{p-1}$.
%}

\begin{figure}
    \centering
    \includegraphics{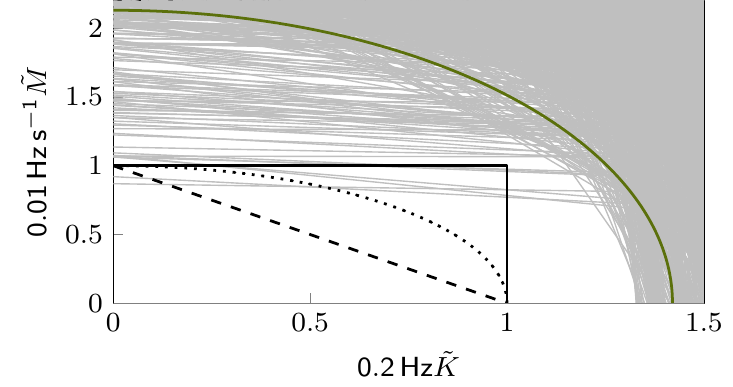}
    \caption{Bounding boxes on reserve provision resulting from Figure~\ref{fig:omegaOmegadot}, for a unit with \SI{1}{\pu}}
    \label{fig:vkvM_constraints}
\end{figure}

Figure \ref{fig:vkvM_constraints} is a map of the $(\omega,\dot{\omega})$ space from Figure~\ref{fig:omegaOmegadot} into $(\vM,\vK)$ space, and with the limits from Figure~\ref{fig:omegaOmegadot} $h$ is $\frac{1}{20}$. Each data point in  $(\omega,\dot{\omega})$ becomes a constraint in $(\vM,\vK)$ and the scaled $\infty$-, 2- and 1-norms translate to scaled 1-, 2- and $\infty$-norms, respectively. Of the three norms, we choose the scaled 1-norm in $(\omega,\dot{\omega})$-space which gives the best fit (like the other norms, it is only violated by few outliers) and which results in a linear and local box-constraint in $(\vM,\vK)$ space, %, the constraint becomes
%\begin{align}
%    &\lVert \cdot \rVert_\infty : &    \begin{bmatrix}\vKi \\ \vMi\end{bmatrix} & \leq \begin{bmatrix} \vPi \\ \vPi \end{bmatrix} \label{eq:vKvM_infnorm} \\
%    &\lVert \cdot \rVert_2 : & \vMi^2 + \vKi^2 & \leq \vPi^2 \label{eq:vKvM_2norm}
%\end{align}
\begin{align}
    \begin{bmatrix}\vK/h \\ \vM\end{bmatrix} & \leq \begin{bmatrix} \vP \\ \vP \end{bmatrix} \label{eq:vKvM_infnorm}
    \eqcomma
\end{align}
which is suitable for efficient linear program formulations.
%Using \eqref{eq:vKvM_infnorm} allows for a linear program formulation and hence is efficiently included in optimizations.%, while \eqref{eq:vKvM_2norm} gives rise to a quadratically constrained program. Also, for given maximum $\omega$ and $\dot{\omega}$, \eqref{eq:vKvM_infnorm} is the least restrictive norm.

\section{Computation of gradient descent directions}\label{sec:derivatives}
% The inertia placement algorithm iteratively searches for parameters  $\vK$, $\vM$ optimizing damping ratio, overshoot and \ac{ROCOF}.
Our inertia placement algorithm (presented formally in Section \ref{sec:optimization}) searches for parameters $\vK$, $\vM$ to optimize damping ratio, overshoot and \ac{ROCOF} subject to constraints.
Hence, we need the gradient or sensitivity of these performance indices with respect to these parameters. The following subsections describe how to compute or approximate these. % overshoot, \ac{ROCOF} and their sensitivities.

\subsection{Computation of the sensitivity of the damping ratio}
The sensitivity of the damping ratio \eqref{eq:dampingratio} with respect~to~$\alpha$~is
%\begin{multline}
%    \D{\damprat}{\alpha} = \D{}{\alpha} \left( \frac{-\reeig}{\sqrt{(\reeig)^2 + (\imeig)^2}} \right) \\
\begin{equation} 
  \D{\damprat}{\alpha}
    = \imeig\frac{\left(\reeig\D{\imeig}{\alpha}-\imeig\D{\reeig}{\alpha}\right)}{\left((\reeig)^2 + (\imeig)^2\right)^{\frac{3}{2}}} \eqcomma
  \label{eq:Da_dzeta}
\end{equation}
%\end{multline}
% To compute \eqref{eq:Da_dzeta} we need the eigenvalue derivatives, which are known to be\cite{Murthy1988}
where the eigenvalue derivatives are obtained from \cite{Murthy1988} as
\begin{align}
    \D{\eig{i}}{\alpha} &= \levecT (\D{\mat{A}}{\alpha}) \revec \eqdot \label{eq:dlambda_dalpha}
\end{align}

\subsection{Computation of the overshoot and its sensitivity}
To find the overshoot \eqref{eq:Mp}, we use the Newton method to search for an extremum of the step response $\mat{y}(\tp)$
\begin{align}
    \tp[\nu+1] &= \tp[\nu] - \left.\left(\D{\mat{y}}{t} \hadad \Dtwo{\mat{y}}{t}\right)\right|_{\vec{t}=\tp[\nu]} \eqcomma \label{eq:NewtonS}
\end{align}
with the n-th derivatives of \eqref{eq:y_of_t} given by
\begin{align}
    \left.\Dn{\vec{y}}{t}\right|_{\vec{t}} & = 
    -\sum\nolimits\nolimits_i (\eig{i})^{(n-1)} \res \! \hadam \eul^{\eig{i} \vec{t}} \eqcomma \label{eq:dy_dt}
    %\D{\vec{y}}{t} &= - \sum\nolimits\nolimits_i \res \eul^{\eig{i} t} \label{eq:dy_dt} \\
	%\Dtwo{\vec{y}}{t} & = -\sum\nolimits\nolimits_i \res \eig{i} \eul^{\eig{i} t} \\
    %\Dthree{\vec{y}}{t} & = -\sum\nolimits\nolimits_i \res (\eig{i})^2 \eul^{\eig{i} t} \eqdot
\end{align}
and where $\hadam$ and $\hadad$ denote the element-wise (Hadamard) multiplication and division.
As $\nu \to \infty$, we obtain the (local) extremum $\mat{y}(\tp)$ at $\tp[\nu] \to \tp$.
%The Newton method converges to a point where $\D{y}{t}$ is zero. 
Since there may be multiple extrema, one needs a good starting point. Gridding the step response and starting from the largest point found leads to the correct extremum if the grid is chosen sufficiently small, e.g., twice the frequency of the highest mode frequency. If an estimate of $\tp$ is available, e.g., from a previous placement iteration, it can be used to initialize the Newton search.
%In following iterations of the inertia placement, the previous peak time can be used as starting point.

It is quite involved to find the sensitivity of the overshoot with respect to $\alpha$: while it is easy to compute the change of the step response with respect to $\alpha$ given a fixed $\mat{t}$, the time of the overshoot $\tp$ is also a function of $\alpha$. Hence, simply taking $\left.\D{y}{\alpha}\right|_{\tp}$ is incorrect. The correct derivative is
\begin{multline}
   \D{\Mp}{\alpha} = \sum\nolimits\nolimits_i \Bigg[\D{\kres}{\alpha} \hadam \left( 1-\eul^{\eig{i}\tp} \right) - \\ \kres \hadam \left( \left(\D{\eig{i}}{\alpha}\right)  \tp + \eig{i} \D{\tp}{\alpha} \right) \hadam \eul^{\eig{i}\tp}\Bigg]
   \label{eq:dalpha_dMp}
\end{multline}
which includes the derivative of the \textit{residues}, the derivative of the \textit{eigenvalues}, and the derivative of the \textit{peak time}.

The \emph{derivative of the residues} $\res$ is given by
\begin{align}
   \D{\res}{\alpha} &= \mat{C} \D{}{\alpha}\left(\revec\levecT\right) \mat{B} \label{eq:dres_dalpha}\\% + \vec{c} \revec\levecT \D{\vec{b}}{\alpha} \\
   \D{\kres}{\alpha} &= \frac{\left(\D{\res}{\alpha}\right) \eig{i} - \res \D{\eig{i}}{\alpha}}{(\eig{i})^2} \eqcomma
\end{align}
%For $\D{\res}{\alpha}$, we need $\D{}{\alpha}\left(\revec\levecT\right)$. This can be computed by \cite{Murthy1988, Borsche2015CDC}
where $\D{}{\alpha}\left(\revec\levecT\right)$  can be computed as follows \cite{Murthy1988, Borsche2015CDC}:
\begin{align}
     \D{}{\alpha}\left(\revec\levecT\right) &= \mkern-9mu \sum\nolimits\nolimits_{j\in \mathcal{N} \setminus i} \mkern-6mu \left[ \revec[j] c_{ij}^\alpha \levecT - \revec c_{ji}^\alpha \levecT[j]\right] \label{eq:duv_dalpha} \\
    c_{ij}^\alpha &= \frac{\levecT[j] (\D{\mat{A}}{\alpha})\revec[i]}{\left(\eig{i}-\eig{j}\right) } \quad i \neq j \eqdot \label{eq:ckj}
\end{align}
Note that for the correct value of $c_{ij}^\alpha$ we used the normalization \eqref{eq:evnormalization}.
Observe that the term $\eig{i}-\eig{j}$ is zero for double eigenvalues. While these usually do not occur in power systems unless the system is perfectly symmetric, one should ensure that the system at hand is well posed.

The \emph{derivative of the peak time} $\tp$ cannot be exactly computed, as $\tp$ is found with the  Newton method. We use as an approximation the derivative of the Newton update \eqref{eq:NewtonS}
\begin{multline}
    \D{\tp}{\alpha} \approx - \D{}{\alpha} \left(\D{\mat{y}}{t} \hadad \Dtwo{\mat{y}}{t}\right) \\
    = - \Big(\left[
            \left(\D{}{\alpha}\D{\mat{y}}{t}\right) \hadam \Dtwo{\mat{y}}{t} - \D{\mat{y}}{t} \hadam \left( \D{}{\alpha} \Dtwo{\mat{y}}{t} \right)
        \right]\hadad
        \\
            % \left( \Dtwo{\mat{y}}{t} \right) ^2
        \left.
            \hadad \left[\Dtwo{\mat{y}}{t} \hadam \Dtwo{\mat{y}}{t}\right]
        \Big)\right|_{\mat{t}=\tp} \label{eq:dtp_dalpha} ,
\end{multline}
for which we need \eqref{eq:dy_dt} and the derivatives of the step response with respect to $\alpha$ explicitly given by:
\begin{align}
    \D{}{\alpha} \D{\mat{y}}{t} &= -\sum\nolimits_i \left[\D{\res}{\alpha} + \res\!\hadam \mat{t} \D{\eig{i}}{\alpha} \right] \hadam \eul^{\eig{i} \mat{t}} \label{eq:DaDt_y}
    \\
    \D{}{\alpha} \Dtwo{\mat{y}}{t} &= -\sum\nolimits_i \left[ \D{\res}{\alpha} \eig{i} + \res \D{\eig{i}}{\alpha} + \right. 
    \nonumber\\ 
    &\left. + \res\!\hadam \mat{t} \eig{i} \D{\eig{i}}{\alpha} \right] \hadam \eul^{\eig{i} \mat{t}} \eqdot \label{eq:DaDt2_y}
\end{align}%
Observe that the approximation \eqref{eq:dtp_dalpha} is actually exact if the Newton method converges in one step. This is the case if $\mat{y}$ is a quadratic function. In our case, $\mat{y}$ consists of sinusoidal functions, and we assume to be at an extremum for the previous parameter value $\alpha$ where sinusoidal functions are described up to fourth order terms by quadratic functions. 
%which is presumably close to the extremum we search for. If we look at the Taylor expansion around an extremum of a sinusoidal function, e.g.,  $\cos (0)$, we see that it is nearly quadratic for small $x$
%\begin{equation}
%     \cos(x) = 1 - \frac{x^2}{2!} + \frac{x^4}{4!} - \frac{x^6}{6!} \cdots \eqdot
%\end{equation}
% While this suggest \eqref{eq:dtp_dalpha} to be a close approximation of the correct derivative for $\tp$ one should keep in mind that 
%we 1)~deal with the sum of many sinusoidal functions, this sum may no longer be close to quadratic, and 2)~we 
%cannot give guarantees on how far we move from our initial optimum by changing $\alpha$. However, 
%While there are no guarantees on the quality of this approximation, it shows 
Hence, in practice, we observe that the approximation \eqref{eq:dtp_dalpha} performs very well; see later simulations in Section~\ref{sec:results}.

\subsection{Computation of the \ac{ROCOF} and its sensitivity}
% The \ac{ROCOF} is the time derivative of $y(t)$, given by \eqref{eq:dy_dt}. 
The \ac{ROCOF} $\ROCOF$ \eqref{eq:Mp} is also found using the Newton method
\begin{align}
    \tROCOF[\nu+1] &= \tROCOF[\nu] - \left.\left(
        \Dtwo{}{t}\mat{y} \hadad
        \Dthree{}{t}\mat{y}
    \right)\right|_{\mat{t}=\tROCOF[\nu]} \eqcomma
\end{align}
where $\tROCOF[\nu] \to \tROCOF$ as $\nu \to \infty$. and $\ROCOF=\mat{y}(\tROCOF)$. 

To find the sensitivity of $\ROCOF$ with respect to $\alpha$ we  use the same approach as for the sensitivity of the overshoot $\D{\Mp}{\alpha}$; see equations \eqref{eq:dalpha_dMp} to \eqref{eq:DaDt2_y}. The RoCoF sensitivity is
\begin{multline}
    \D{\ROCOF}{\alpha} = - \sum\nolimits_i \Bigg[ \D{\res}{\alpha} + \res\!\hadam \tROCOF \D{\eig{i}}{\alpha} + \\ 
    + \res\!\hadam \left(\D{\tROCOF}{\alpha}\right) \eig{i} \Bigg] \hadam \eul^{\eig{i}\tROCOF} \eqcomma
\end{multline}
using \eqref{eq:dlambda_dalpha} and \eqref{eq:dres_dalpha}. The sensitivity of $\tROCOF$ is approximated by the derivative of the \ac{ROCOF} Newton update step,
\begin{multline}
    \D{\tROCOF}{\alpha} \approx - \D{}{\alpha} \left( \Dtwo{\mat{y}}{t} \hadad \Dthree{\mat{y}}{t} \right)\\
    = - \Big(\left[
        \left(\D{}{\alpha}\Dtwo{\mat{y}}{t}\right) \hadam \Dthree{\mat{y}}{t} - \Dtwo{\mat{y}}{t} \hadam \left( \D{}{\alpha} \Dthree{\mat{y}}{t} \right)\right]\hadad
        \\ 
        \left.\hadad \left[ \Dthree{\mat{y}}{t} \hadam \Dthree{\mat{y}}{t} 
    \right]\Big)\right|_{\mat{t}=\tROCOF} ,
\end{multline}
with \eqref{eq:dy_dt}, \eqref{eq:DaDt2_y} and 
\begin{multline}
    \D{}{\alpha} \Dthree{\mat{y}}{t} = -\sum\nolimits_i \Big[ \D{\res}{\alpha} \left(\eig{i}\right)^2 
    + 2 \res \eig{i} \D{\eig{i}}{\alpha} + \\ + \res\!\hadam\mat{t} \left(\eig{i}\right)^2 \D{\eig{i}}{\alpha} \Big] \hadam \eul^{\eig{i} \mat{t}} . 
\end{multline}

\section{Optimal inertia placement algorithm}\label{sec:optimization}
In the  following we present the objectives of our synthetic inertia allocation optimization (Section \ref{ssec:objective}), pose the program in a general form (Section \ref{ssec:algorithm}), and a sequential linear programming approach to solve it (Section~\ref{ssec:iteration}).

\subsection{Formulation of optimization objectives}\label{ssec:objective}
We aim to co-optimize three metrics, namely the damping ratio $\zeta$, the overshoot $\Mp$ and the \ac{ROCOF} $\ROCOF$. 

To maximize the smallest $\zeta$, we introduce the variable $\zetamin$
\begin{subequations}
\begin{align}
    \zetamin     & \leq \damprat && \forall \{i|\imeig[i] > 0\} \eqcomma
\end{align}
and consider the cost term
\begin{align}
    \min\nolimits_\alpha -\czeta \zetamin 
\end{align}
\end{subequations}
with positive $\czeta$, pushing $\zetamin$ against the smallest $\damprat$.

Similarly, to minimize the steepest \ac{ROCOF}, $\ROCOFmax$, we consider the cost term
\begin{subequations}
\begin{equation}
    \min\nolimits_\alpha \cROCOF \ROCOFmax 
\end{equation}
with positive parameter $\cROCOF$ and subject to
\begin{align}
    \ROCOFmax &\geq  |\ROCOFij| && \forall i,j~.
\end{align}
\end{subequations}
Analogous constraints and costs are used for the overshoot~$\Mp$.

\subsection{Inertia and damping placement algorithm}\label{ssec:algorithm}
We pose the optimal synthetic inertia and damping placement problem as a multi-objective optimization problem:
\begin{subequations}
\label{eq:opt}
\begin{align}
    \min\nolimits_\alpha -\czeta \zetamin +\cROCOF \ROCOFmax + \cMp \Mpmax \label{eq:opt_cost}
\end{align}
\begin{align}
    \text{s.t.} 
    && \underline{\zeta}  &\leq \damprat                        & \forall \{i|\imeig[i] > 0\} \label{eq:opt_zetabound}\\
    && \underline{\ROCOF} &\leq \ROCOFij \leq \overline{\ROCOF} & \forall i,j \label{eq:opt_Rbound}\\
    && \underline{\Mp}    &\leq \Mpij \leq \overline{\Mp}       & \forall i,j \label{eq:opt_Sbound}\\[6pt]
    && 0                  &\leq \vK_v/h \leq \vPi                 & \forall v \label{eq:opt_Kbound}\\
    && 0                  &\leq \vM_v \leq \vPi                 & \forall v \label{eq:opt_Mbound}
\end{align}%
\end{subequations}%
The cost function \eqref{eq:opt_cost} combines the three-level objective discussed in Section~\ref{ssec:objective}. The constraints \eqref{eq:opt_zetabound} to \eqref{eq:opt_Sbound} put strict bounds on all damping ratios, \ac{ROCOF}s and overshoots. Constraints \eqref{eq:opt_Kbound} and \eqref{eq:opt_Mbound} define bounds on synthetic inertia and damping provision according to Section~\ref{sec:vMvK_constraints}; confer \eqref{eq:vKvM_infnorm}.

Thus, our inertia and damping allocation problem \eqref{eq:opt} combines system-level objectives, such as the damping ratio, together with explicit time-domain criteria and strict constraints on \ac{ROCOF}, overshoot, and device-level power limits.

\subsection{Sequential linear programming approach}\label{ssec:iteration}
The optimization problem \eqref{eq:opt} is non-linear, typically large-scale for the considered system, and highly non-convex. Hence we use a sequential linear programming approach iterating over parameters $\alpha$ until we reach a local optimum. 

At each iteration a first-order (linear) approximation of \eqref{eq:opt} is obtained as follows. Given values $\alpha^\nu$, $\zetan[j]$, $ \ROCOF[\nu]$, and $\Mp[\nu]$ at iteration $\nu$, the performance metrics are then updated by means of the sensitivities derived in Section~\ref{sec:derivatives} as
\begin{subequations}%
\label{eq:n1}%
\begin{align}%
    \zetanone[j] & = \zetan[j] + \sum\nolimits_{\forall \alpha} \D{\zetan[j]}{\alpha} \Delta \alpha^\nu  && \forall \{j|\imeig[j] > 0\} \label{eq:zetan1}\\
    % \ROCOFi[\nu+1] & = \ROCOFi[\nu] + \sum\nolimits_{\forall \alpha} \D{\ROCOFi[\nu]}{\alpha} \Delta \alpha^\nu && \forall d,b \\
    % \Mpi[\nu+1] & =\Mpi[\nu] + \sum\nolimits_{\forall \alpha} \D{\Mpi[\nu]}{\alpha} \Delta \alpha^\nu && \forall d,b \label{cs:Pn1}
    \ROCOFnone[\nu+1] & = \ROCOF[\nu] + \sum\nolimits_{\forall \alpha} \D{\ROCOF[\nu]}{\alpha} \Delta \alpha^\nu \label{eq:Rn1} \\
    \Mpnone[\nu+1] & =\Mp[\nu] + \sum\nolimits_{\forall \alpha} \D{\Mp[\nu]}{\alpha} \Delta \alpha^\nu \label{eq:Sn1} \eqcomma 
\end{align}%
\end{subequations}%
with $\alpha^\nu + \Delta \alpha^\nu = \alpha^{\nu + 1}$, $\Delta \alpha^\nu$ being a (time-varying) step-size, and $\zetan[j]=\zeta(\alpha^\nu)$, $\ROCOF[\nu] = \ROCOF(\alpha^\nu)$ and $\Mp[\nu] = \Mp(\alpha^\nu)$. 
%The left-hand terms in \eqref{eq:n1} are linear approximations, and by letting $\zetanone[j] = \damprat[j]$, $\ROCOFnone[\nu+1] = \ROCOF$ and $\Mpnone[\nu+1] = \Mp$ in \eqref{eq:opt}, a linear program is obtained.

The left-hand terms in \eqref{eq:n1} are first-order approximations when $\alpha^{\nu}$ is updated to $\alpha^{\nu + 1}$. By setting $\zetanone[j] = \damprat[j]$, $\ROCOFnone[\nu+1] = \ROCOF$ and $\Mpnone[\nu+1] = \Mp$ in \eqref{eq:opt}, we obtain a linear programming formulation. % is obtained.

\textit{Updates and limits of the step size:}
As \eqref{eq:n1} are only locally valid linearizations, we need to limit the step size $\Delta \alpha^{\nu}$ as
\begin{align}
	-\Delta \alpha^\mathrm{max} & \leq  \Delta \alpha^\nu \leq \Delta \alpha^\mathrm{max} && \forall \alpha
	\eqcomma
\end{align}%
where $\Delta \alpha^\mathrm{max}>0$.
After each iteration, the updated system matrix $\mat{A}(\alpha^{\nu+1})$ and performance indices $\zeta$, $\Mp$ and $\ROCOF$ are computed. If they show an improvement, $\alpha^{\nu+1}$ is kept. Otherwise, the previous value $\alpha^{\nu}$ is used, and the step size $\Delta \alpha^\mathrm{max}$ is halved for all $\Delta \alpha$ that hit $\Delta \alpha^\mathrm{max}$.

\indent\textit{Iterations:}
Due to the mismatch between the linearly approximation $\ROCOFnone[\nu+1]$ and the true value $\ROCOF(\alpha^{\nu+1})$, the new starting point $\ROCOF(\alpha^{\nu+1})$ may violate the constraint \eqref{eq:opt_Rbound}. To ensure feasibility, we add a slack variable to this constraint,
\begin{subequations}
\begin{align}
	\underline{\ROCOF} & \leq \ROCOFij[\nu+1] - \epsROCOF \leq \overline{\ROCOF} && \forall i,j
\end{align}
The slack variable is only added if the starting point of the iteration is infeasible
\begin{align}
	0 & \leq \epsROCOF && \forall \{i,j|\ROCOFij[\nu]> \overline{\ROCOF}\} \\
    0 & =    \epsROCOF && \forall \{i,j|\ROCOFij[\nu] \in [\underline{\ROCOF}, \overline{\ROCOF}]\}\\
	0 & \geq \epsROCOF && \forall \{i,j|\ROCOFij[\nu] < \underline{\ROCOF}\}
\end{align}
and the slack is penalized with a large cost term $\cepsROCOF >0$ as
\begin{equation}
    \min\nolimits_\alpha \sum\nolimits_{d,b} \cepsROCOF \epsROCOF \eqdot
\end{equation}
\end{subequations}
The same approach is used to ensure feasibility of $\zeta$ and $\Mp$.

\textit{Stopping criterion:}
The algorithm terminates after a fixed number of iterations, when the performance improvement is smaller than a threshold, or when the step size for all $\alpha$ is below a threshold.

\subsection{Considerations on numerics}
The main computational effort at each iteration of our algorithm is to obtain all eigenvalues, eigenvectors and sensitivities. Clearly, this is computationally burdensome for large systems, but the scaling is reasonable as shown below.

Of the needed computations, the derivatives of the eigenvectors have the worst scaling. For each $\alpha$, we need to compute all $N$ eigenvector derivatives. For each eigenvector derivative, we need $c_{kj}^\alpha$ which has dimension $\Rdim{N \times N}$. Assuming $N^A$ parameters, this leads to $N_A \, N^3$ entries. By exploiting the structure of our problem, we can significantly reduce the number of entries: while we need all eigenvector derivatives, we only need the entries that correspond to non-zero entries in $\mat{B}$ and $\mat{C}$. These scale with the number of disturbances $N_B$ and observed frequencies $N_C$, giving a scaling of $\mathcal{O} (N \, N_A \, N_B \, N_C) $. The dimensions $N_B$ and $N_C$ are much smaller than $N$. Additionally, if we change the model and add some states while keeping the number of disturbances and outputs constant, the computational effort for the eigenvector derivatives scales linearly instead of cubic. This allows much more freedom on modeling choices.

% With $N$ the dimension of the system, $A$ the number of parameters, $D$ the number of disturbances and $B$ the number of observed bus frequencies, the residuals $\res$ have dimension $\Rdim{N\times D \times B}$, and their derivatives $\D{\res}{\alpha}$ are in $\Rdim{A\times N \times D \times B}$. The most expensive and worst scaling operation are the derivatives of the eigenvectors. For each $\alpha$, we need to compute all $N$ eigenvector derivatives. For each eigenvector derivative, we in turn need $c_{kj}^\alpha$ which has dimension $\Rdim{N \times N}$. In total, this leads to $A \, N^3$ entries. Using the structure of our problem, we can significantly reduce the number of entries: while we need all eigenvector derivatives, we only need the entries that correspond to non-zero entries in $\mat{B}$ and $\mat{C}$. These scale with the number of disturbances $D$ and observed frequencies $B$, giving a scaling of $\mathcal{O} (N \, A \, D \, B) $. $D$ and $B$ are not only much smaller than $N$. If we change the model and add some states, the computational effort for the eigenvector derivatives will scale linearly instead of cubic. This allows much more freedom on modeling choices.

Computation of the overshoot $\Mp$ and \ac{ROCOF} $\ROCOF$ matrices, both in $\Rdim{D\times B}$, via the iterative Newton method, is comparably efficient. While at the first iteration of the overall placement algorithm, we grid the system to find the global extrema, we use the previous $\tp$ and $\tROCOF$ as starting points for the Newton method in the next iteration and compute all $\tp$ and $\tROCOF$ in parallel. With this implementation, it takes usually only few Newton iterations to find the  values of $\Mp$ and $\ROCOF$.

%\subsection{Alternative cost on $\zeta$, $\Mp$ or $\ROCOF$}
\subsection{Alternative optimization problem formulations}
\label{subsec: alternatives}

In the following, we briefly discuss a few alternative formulations of the optimal inertia and damping allocation problem \eqref{eq:opt} which are useful in other scenarios.

\textit{Power capacity:}
Previously, we assumed the power capacity $\vPi$ of each inertia device to be fixed. This renders the program \eqref{eq:opt} a scheduling problem answering how much of the available inertia at a certain node should be used, depending on the expected system state. Alternatively, we can make $\vPi$ a decision variable, rendering the program \eqref{eq:opt} a planning problem. If a fixed amount of inertia-devices is to be placed in the system, one would add the budget constraint
\begin{align}
    \sum\nolimits \vPi & \leq \Pbdg \eqdot \label{eq:budget}
\end{align}
Finally, $\Pbdg$ can be itself be a decision variable as well with an associated cost that reflects the investment cost of synthetic inertia devices in a planning program.

\textit{Average performance:}
Instead of optimizing and limiting the worst-case performance as in \eqref{eq:opt}, we could also optimize the average performance objective
\begin{equation}
    \min\limits_{\alpha^{\nu+1}} -c^\zeta \sum\limits_j \zetan[j] + c^\mathrm{S1} \sum\limits_{i,j} \Mpij  + \cROCOFmean \sum\limits_{i,j} \ROCOFij 
    \label{eq: avg perf}
    \eqdot
\end{equation}
Such an average objective allows to trade off damping, \ac{ROCOF} and overshoot between different buses in the system.

\begin{figure}
    \centering
    \includegraphics[scale=0.5]{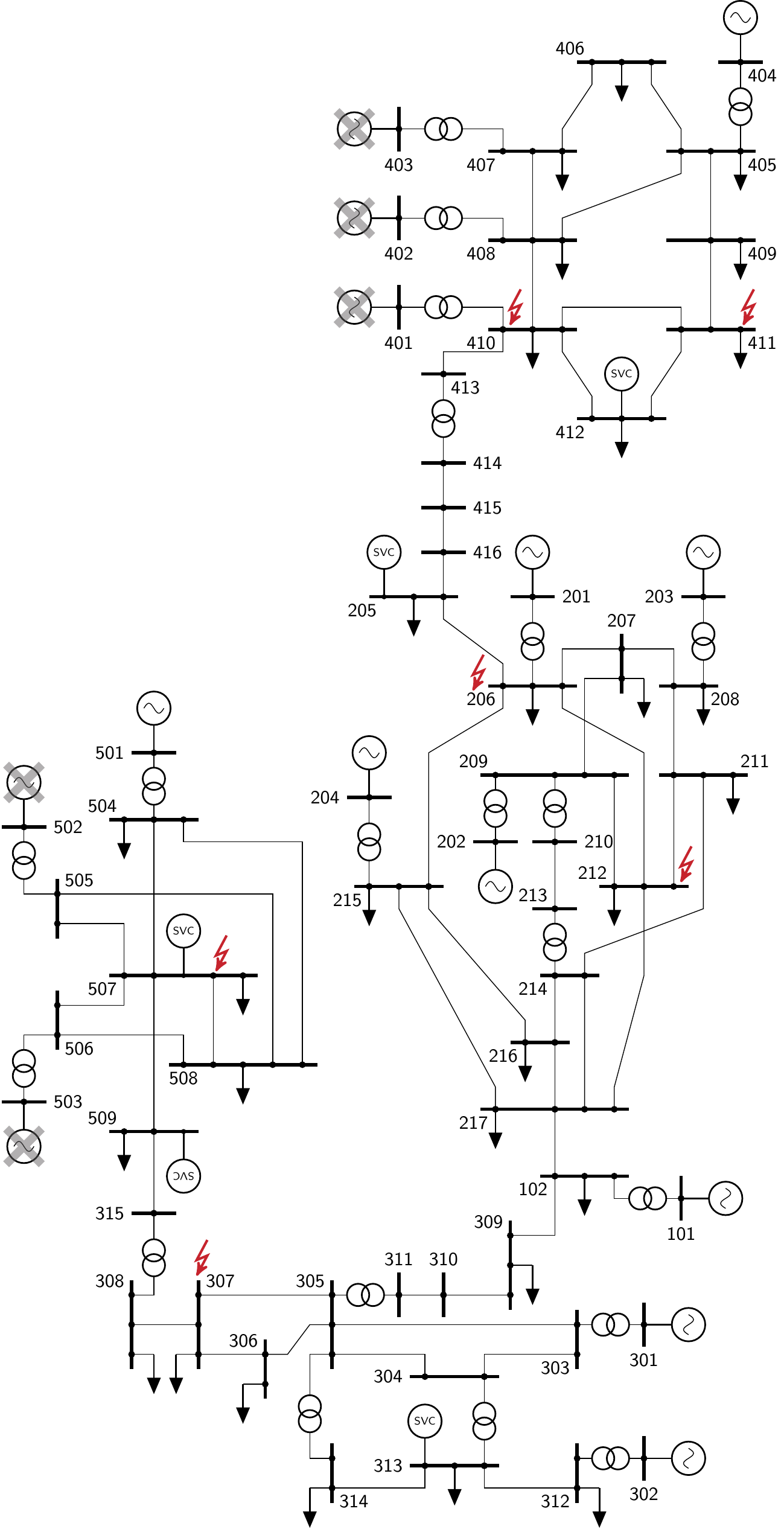}
    \caption{Topology of the test system.}
    \label{fig:topo}
\end{figure}

%\begin{figure}[htbp]
%    \centering
%    \includegraphics{lambda.pdf}
%    \caption{Dominating eigenvalues of the test system with added motor loads. The generator modes (green/blue) and motor modes (black) are highlighted.}
%    \label{fig:lambda0}
%\end{figure}

\begin{figure*}[htbp]
    \centering
    \includegraphics[scale=1]{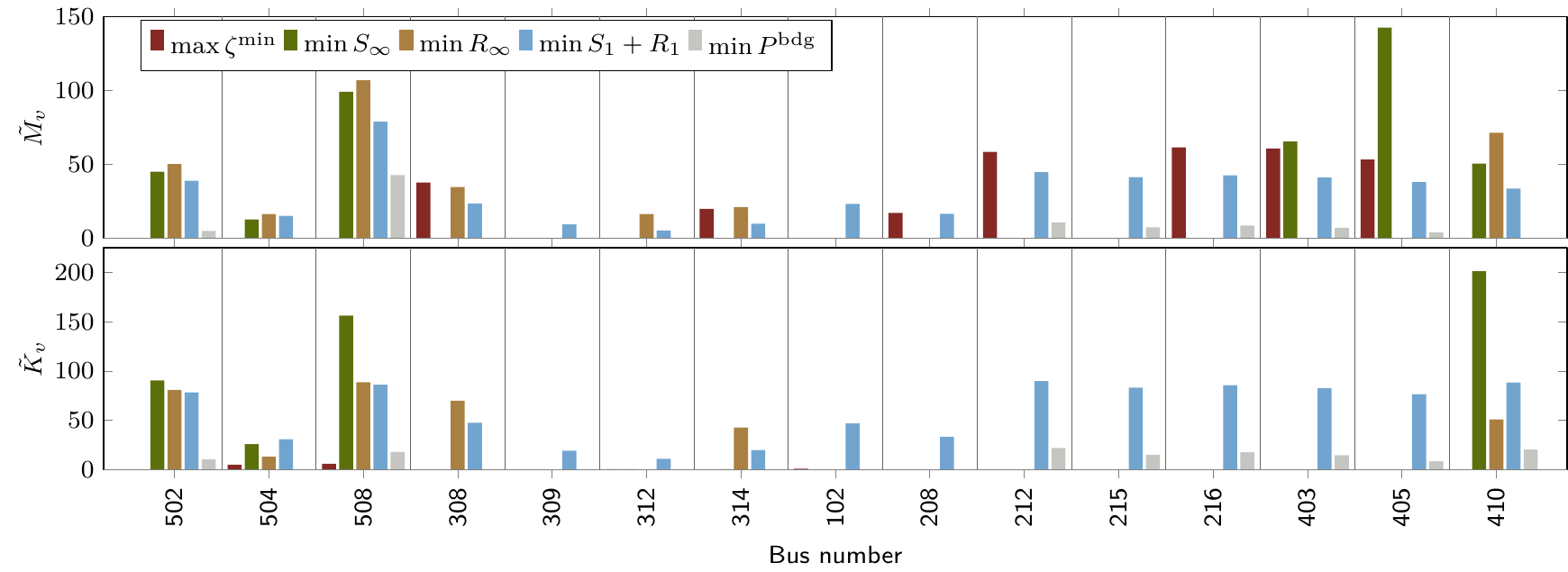}
    \caption{Synthetic inertia and damping distribution with different cost functions}
    \label{fig:distribution}
\end{figure*}

\section{Description of the test case}\label{sec:testsystem}\label{sec:testcase}
We use a modified version of the Australian 14 generator system as a test case to illustrate the utility of our synthetic inertia and damping placement algorithm. %The first subsection introduces the system, the second gives results.

%\subsection{System description}

\textit{Base system:}
The Australian system \cite{Gibbard2014b}, see Figure~\ref{fig:topo}, consists of 14 generators and 59 buses. Gibbard and Vowles describe six load cases, of which we chose the heavy loading case. The system consists of five areas connected in a string-like layout. The main demand, Melbourne and Sydney, is in the middle of the system, while significant generation is located in all areas including the two far ends of the system. We use the AVR and PSS parameters given in \cite{Gibbard2014b}. 

The model is extended with motor loads and load damping. We assume \SI{10}{\percent} of the load to be from motors with an inertia of \SI{1.5}{\second}. Load damping is set to \SI{2.5}{\pu}, and dynamic loads are behind an \SI{0.1}{\pu} inductance \cite{Kundur1994}. 

% Figure~\ref{fig:lambda0} shows the initial eigenvalue distribution, with the modes associated with motors in black, inter-area oscillations in blue and loca area osciallations in green.

\textit{Low-inertia case study:}
For the case study, we remove five generators from the system, namely 401, 402, 403, 502 and 503. These are generators in the West and North ends of the system, where there is abundant wind  and solar resources, respectively, and which are likely areas for RES deployment in Australia. We assume renewable generation to have a grid-following maximum power-point tracking control feeding constant active and reactive power into the system.

\textit{Disturbances:} We model disturbances as sudden load increases of \SI{250}{\mega\watt} at some load buses, namely 206, 212, 307, 410, 411 and 508. We chose such generic faults as they do not affect the $\mat{A}$-matrix of the system. 
% If a generator or line tripping is to be modeled, one would have to consider the effect on the dynamics matrix---which would be a straight-forward but tedious extension to our framework.

\textit{Monitored frequencies:} To identify the effect of removing generators and adding synthetic inertia, we monitor the frequency $\omega^\mathrm{G}$ at all remaining conventional generators and compute $\zeta$, $\ROCOF$ and $\Mp$ at these buses.

%\fdmargin{Do you mean $P^{budget}$ instead of $\Mbdg$ ?}
\textit{Synthetic inertia and damping budget:} For better comparability with the initial system, we allow the same amount of inertia to be added to the system as is lost due to generator removal. In system base this amounts to an inertia budget $\Mbdg$ of \SI{475}{\second}, which depending on the largest (expected) ROCOF translates to the power budget
\begin{equation}
    \Pbdg = \Mbdg \max |\dot{\omega}| \eqdot
\end{equation}
Finally, the choice of $h$ in \eqref{eq:opt_Kbound} is a relevant design parameter, which traces back to the observed system dynamics in \eqref{eq:omegapnorm}. The observed frequency and ROCOF in the CE system suggested $h = \frac{1}{20}$, typical settings for protection relays are at \SI{2}{\hertz} and \SI{0.5}{\hertz\per\second}, suggesting $h=\frac 1 4$, while simulations of the test system give higher ROCOF then frequency excursions, suggesting $h \approx 3$. We have chosen $h=1$ but recommend to assess this carefully for the system under consideration.

\section{Results}\label{sec:results}
In this section we compare different cost functions for our placement algorithm. We test five approaches: 1) maximizing worst-case damping ratio, 2) minimizing worst-case \ac{ROCOF}, 3) minimizing worst-case overshoot, 4) co-optimizing {average} overshoot and \ac{ROCOF} as in \eqref{eq: avg perf}, and 5) penalizing the expenditure of synthetic inertia and damping; see Section~\ref{subsec: alternatives}.

We first discuss and compare the first four case studies.
Figure~\ref{fig:distribution} shows the allocation of synthetic inertia depending on each of the {four} cost functions. It is immediately evident that the different cost functions lead to significantly different inertia distributions. 
Table~\ref{tab:results} gives a comparison and cross-validation of the results for {these four} cases. As performance indices we use the worst-case metrics $\zetamin$, $\ROCOFmax$ and $\Mpmax$; the total allocation of inertia and of damping, $\sum\nolimits \vM$ and $\sum\nolimits \vK$; and the mean \ac{ROCOF} and mean overshoot, $\ROCOFone$ and $\Mpone$.
The optimal placement of inertia outperforms the initial allocation and helps to alleviate the loss of generators. Depending on the cost function, the performance metrics are affected in quite different ways. Also, the inertia budget is never fully utilized, hinting at the fact that with optimal placement, actually little synthetic inertia is needed.
We also observe that each performance metric is lowest when it is considered in the cost function, suggesting that the results are plausible.

\begin{table}[b]
    \centering
    \caption{Comparison of results}\label{tab:results}
    \sisetup{table-format = 3.0, table-auto-round}
    \begin{tabular}{l@{~}S[table-format = 2.1]@{~}S@{~}S[table-format = 2.1]@{~}S@{~}S@{~}S[table-format = 2.0]@{~}S[table-format = 2.0]}
        \toprule
        Metric &  {$\zetamin$} & {$\ROCOFmax$} & {$\Mpmax$} & {$\sum\nolimits \vM$} & {$\sum\nolimits \vK$} & {$\ROCOFone$} & {$\Mpone$}\\
         &  {\%} & {[mHz/s]} & {[mHz]} & {[pu]} & {[pu]} & {[mHz/s]} & {[mHz]} \\\midrule
        initial system          & 18.6414 & 193.108 & 56.4528 & {--} & {--} & 30.2605 & 20.0082 \\
        low inertia             & 19.1432 & 395.756 & 98.2757 & {--} & {--} & 34.9901 & 18.9289 \\ \midrule
        $\min\nolimits -\zetamin$        & 19.2386 & 378.861 & 90.1813 & 310.717 & 6.44673 & 30.7658 & 18.1864 \\
        $\min\nolimits \Mpmax$           & 14.9994 & 96.2151 & 27.3695 & 417.48 & 237.665 & 21.4256 & 12.0297 \\
        $\min\nolimits \ROCOFmax$        & 14.9999 & 94.0839 & 28.4442 & 429.313 & 163.326 & 22.8881 & 15.5199 \\
        $\min\nolimits \Mpone + \ROCOFone$     & 17.0341 & 96.247 & 27.5439 & 464.629 & 438.961 & 18.8276 & 10.7923 \\ \midrule
        $\min\nolimits \sum\nolimits \vPi$ & 15.9 & 95.9 & 30.0 &  87.6 &  63.3 & 24.5 & 15.7 \\
        \bottomrule
    \end{tabular}
    \label{tab:my_label}
\end{table}

\begin{figure}[htbp]
    \centering
    \includegraphics[scale=0.95]{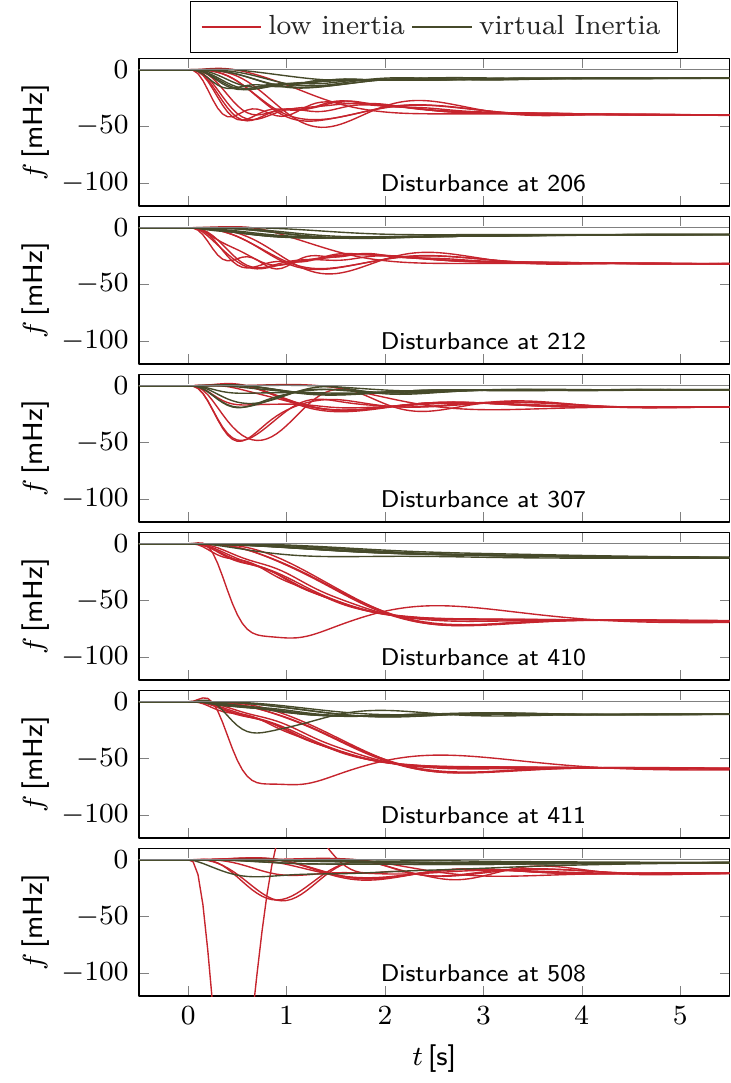}
    \caption{Step response of linearized system. Low inertia (red) and optimal inertia placement (green) with cost function 4).}
    \label{fig:stepResponse}
\end{figure}

\textit{1) Optimizing the damping ratio $\zeta^{\max}$:}
To maximize the damping ratio, our allocation algorithm places inertia mainly in the center of the system and places very little damping. 
% While this counter-intuitive at first, it seems that the time delay in the synthetic inertia response is beneficial for damping inter-area oscillations.

\textit{2) Minimizing the largest overshoot $\Mpmax$:}
While the largest overshoot occurs at bus 501, to minimize $\Mpmax$ inertia and damping are placed in Area~4, and only some in Area~5. It seems that already little additional virtual inertia in Area~5 suffices to alleviate overshoot issues.

\textit{3) Minimizing the largest \ac{ROCOF} $\ROCOFmax$:}
The largest \ac{ROCOF} is found at bus 501, after a disturbance at bus 508. Accordingly, to minimize the \ac{ROCOF} most inertia as well as some damping are allocated in Area~5. Incidentally, this is the area of a recent black-out in the grid, blamed on insufficient fault-ride-through capabilities of wind generation \cite{AEMO2016}.

\textit{4) Co-optimizing the overshoot and the \ac{ROCOF}:}
Penalizing the average overshoot and \ac{ROCOF} {as in \eqref{eq: avg perf}} leads to a very even distribution of damping and inertia at the outer areas 2, 4 and 5. It also leads to $\ROCOFmax$ and $\Mpmax$ that are close to the ones achieved in the two previous approaches. It seems that such a cost function taking into account all frequencies and not only the worst excursions gives the most benign system behaviour. Figure~\ref{fig:stepResponse} shows the step response with this distribution compared to the low-inertia case.

\textit{5) Minimal inertia expenditure:}
Finally, we minimize the use of synthetic inertia and damping, {by making $\Pbdg$ in \eqref{eq:budget} a decision variable,} while keeping the \ac{ROCOF}, overshoot and damping ratio in bounds of \SI{100}{\milli\hertz\per\second}, \SI{30}{\milli\hertz} and \SI{15}{\percent}. {Perhaps surprisingly, we can reduced synthetic inertia requirements by a factor five within the imposed constraints.}
%, at the cost of a slightly higher \ac{ROCOF} and overshoot.

%\begin{table}[ht]
%    \sisetup{table-format = 3.1, table-auto-round}
%    \centering
%    \caption{}
%    \label{}
%    \begin{tabular}{l@{\,}lSSSSSS}
%        \toprule
%        & & & & {$ -\zetamin$} & {$ \Mpmax$} & {$ \ROCOFmax$}  & {$ \Mp + \ROCOF$} \\
%        \midrule
%        $\zetamin$ & \% & 18.6414 & 19.1432 & 19.2386 & 14.9994 & 14.9999 & 17.0341 \\
%        $\ROCOFmax$ & [mHz/s] & 193.108 & 395.756 & 378.861 & 96.2151 & 94.0839 & 96.247 \\
%        $\Mpmax$ & [mHz] & 56.4528 & 98.2757 & 90.1813 & 27.3695 & 28.4442 & 27.5439 \\
%        $\sum\nolimits \vM$ & [pu] & {--} & {--} & 310.717 & 417.48 & 429.313 & 464.629 \\
%        $\sum\nolimits \vK$ & [pu] & {--} & {--} & 6.44673 & 237.665 & 163.326 & 438.961 \\
%        $\sum\nolimits \ROCOF$ & [mHz/s] & 30.2605 & 34.9901 & 30.7658 & 21.4256 & 22.8881 & 18.8276 \\
%        $\sum\nolimits \Mp$ & [mHz] & 20.0082 & 18.9289 & 18.1864 & 12.0297 & 15.5199 & 10.7923 \\
%        \bottomrule
%    \end{tabular}
%\end{table}

%\fdmargin{I would omit this sentence unless they ask for it. Keep some easy targets for the reviewers}
\section{Conclusion}
This paper presented an algorithm for optimal inertia placement with explicit time-domain constraints. A case study on the Australian grid shows the applicability to realistic power system models. 
%Since the optimization is based on a linearized model, it is necessary to check the results against a non-linear model of the system.
With that in mind, the algorithm can be a valuable tool in short and long term planning of power system stability and inertia deployment.

The approach can be easily extended to answer more detailed questions. For example, one can extend the modeling framework to include HVDC lines that emulate inertia by transferring energy from one part of the system to another. 
%In this case, one would get have to change \eqref{eq:vM_feedback} to something of the form
%\begin{equation}
%    \mat{\tilde{C}_v} = \begin{bmatrix} 1 & 0 \\ 0 & -1 \end{bmatrix}
%\end{equation}
%to describe the energy transfer, but the rest of the algorithm would stay unchanged.

%\fdmargin{same: I would avoid this and keep it in case they ask for it}
%We avoided modeling disconnection of generators, as that would lead to a changed system matrix after such a disturbance. Apart from a tedious implementation, it is straightforward to adjust the algorithm accordingly. The same is true if one is interested in faults which would lead to a system split into two or more areas.

Another direction of research is valuation of inertia provision. The optimization problem gives rise to a notion of location marginal inertia prices in line with traditional locational marginal pricing, as argued in \cite{Borsche2016Diss}.

\section*{Acknowledgments}
The authors would like to thank G\"oran Andersson for the continued support, Luis Ruoco for discussions on the small signal modeling and provision of the small signal toolbox, Micheal Gibbard and David Vowles for providing the original version of the Australian test system, and Tao Liu and David Hill for the collaboration that started this research.

\bibliographystyle{IEEEtran}
\bibliography{library,new}

\begin{thebibliography}{10}
\providecommand{\url}[1]{#1}
\csname url@samestyle\endcsname
\providecommand{\newblock}{\relax}
\providecommand{\bibinfo}[2]{#2}
\providecommand{\BIBentrySTDinterwordspacing}{\spaceskip=0pt\relax}
\providecommand{\BIBentryALTinterwordstretchfactor}{4}
\providecommand{\BIBentryALTinterwordspacing}{\spaceskip=\fontdimen2\font plus
\BIBentryALTinterwordstretchfactor\fontdimen3\font minus
  \fontdimen4\font\relax}
\providecommand{\BIBforeignlanguage}[2]{{%
\expandafter\ifx\csname l@#1\endcsname\relax
\typeout{** WARNING: IEEEtran.bst: No hyphenation pattern has been}%
\typeout{** loaded for the language `#1'. Using the pattern for}%
\typeout{** the default language instead.}%
\else
\language=\csname l@#1\endcsname
\fi
#2}}
\providecommand{\BIBdecl}{\relax}
\BIBdecl

\bibitem{Ulbig2014IFAC}
A.~Ulbig, T.~S. Borsche, and G.~Andersson, ``{Impact of Low Rotational Inertia
  on Power System Stability and Operation},'' in \emph{Proceedings of the 19th
  IFAC World Congress}, Cape Town, aug 2014, pp. 7290--7297.

\bibitem{tielens2016relevance}
P.~Tielens and D.~Van~Hertem, ``The relevance of inertia in power systems,''
  \emph{Renewable and Sustainable Energy Reviews}, vol.~55, pp. 999--1009,
  2016.

\bibitem{EirGrid2012}
EirGrid and Soni, ``{DS3: System Services Review TSO Recommendations},''
  EirGrid, Tech. Rep., 2012.

\bibitem{Ercot2013}
ERCOT, ``{Future Ancillary Services in ERCOT},'' ERCOT, Tech. Rep., 2013.

\bibitem{statnett2016inertia}
``Challenges and opportunities for the nordic power system,'' statnett,
  fingrid, energinet.dk, svenska kraftn\"att, Tech. Rep., Aug 2016.

\bibitem{Bevrani2014}
H.~Bevrani, T.~Ise, and Y.~Miura, ``Virtual synchronous generators: A survey
  and new perspectives,'' \emph{Intl. Journal of Electrical Power \& Energy
  Systems}, vol.~54, Jan. 2014.

\bibitem{DG-SB-BKP-FD:17}
D.~Gross, S.~Bolognani, B.~K. Poolla, and F.~D{\"o}rfler, ``Increasing the
  resilience of low-inertia power systems by virtual inertia and damping,'' in
  \emph{Bulk Power Systems Dynamics and Control Symposium (IREP)}, 2017, to
  appear.

\bibitem{Rakhshani2016}
E.~Rakhshani, D.~Remon, A.~M. Cantarellas, and P.~Rodriguez, ``Analysis of
  derivative control based virtual inertia in multi-area high-voltage direct
  current interconnected power systems,'' \emph{IET Generation, Transmission \&
  Distribution}, vol.~10, no.~6, pp. 1458--1469, 2016.

\bibitem{poolla2016placing}
B.~K. Poolla, S.~Bolognani, and F.~D{\"o}rfler, ``Placing rotational inertia in
  power grids,'' in \emph{American Control Conference (ACC), 2016}.\hskip 1em
  plus 0.5em minus 0.4em\relax IEEE, 2016, pp. 2314--2320.

\bibitem{Pirani2017}
M.~Pirani, J.~W. Simpson-Porco, and B.~Fidan, ``System-theoretic performance
  metrics for low-inertia stability of power networks,'' \emph{arXiv preprint
  arXiv:1703.02646}, 2017.

\bibitem{Mesanovic2016}
A.~Mesanovic, U.~M{\"{u}}nz, and C.~Hyde, ``{Comparison of H∞ , H2 , and pole
  optimization for power system oscillation damping with remote renewable
  generation},'' in \emph{IFAC Workshop on Control of Transmission and
  Distribution Smart Grids - CTDSG’16}, Prague, 2016.

\bibitem{Borsche2015CDC}
T.~S. Borsche, T.~Liu, and D.~J. Hill, ``{Effects of Rotational Inertia on
  Power System Damping and Frequency Transients},'' in \emph{54th IEEE
  Conference on Decision and Control (CDC)}, Osaka, 2015.

\bibitem{Vournas1987}
C.~D. Vournas and B.~C. Papadias, ``{Power system stabilization via parameter
  optimization-application to the Hellenic interconnected system},'' \emph{IEEE
  Transactions on Power Systems}, vol.~2, no.~3, pp. 615--622, 1987.

\bibitem{BKP-DG-TB-SB-FD:17}
B.~K. Poolla, D.~Gross, T.~Borsche, S.~Bolognani, and F.~D{\"o}rfler, ``Virtual
  inertia placement in electric power grids,'' in \emph{Energy Markets and
  Responsive Grids}, J.~Stoustrup, Ed., 2017.

\bibitem{Gibbard2014b}
M.~Gibbard and D.~Vowles, ``{Simplified 14-Generator Model of the South East
  Australian Power System, Revision 4},'' Tech. Rep. June, 2014.

\bibitem{AEMO2016}
AEMO, ``{Update Report - Black System Event in South Australia on 28 September
  2016},'' Tech. Rep., 2016.

\bibitem{Bergen1981}
A.~R. Bergen and D.~J. Hill, ``{A structure preserving model for power system
  stability analysis},'' \emph{IEEE Transactions on Power Apparatus and
  Systems}, vol. PAS-100, no.~1, pp. 25--35, 1981.

\bibitem{FD-FB:11d}
F.~D{\"o}rfler and F.~Bullo, ``{K}ron reduction of graphs with applications to
  electrical networks,'' \emph{IEEE Transactions on Circuits and Systems I:
  Regular Papers}, vol.~60, no.~1, pp. 150--163, January 2013.

\bibitem{Murthy1988}
D.~V. Murthy and R.~T. Haftka, ``{Derivatives of eigenvalues and eigenvectors
  of a general complex matrix},'' \emph{International Journal for Numerical
  Methods in Engineering}, vol.~26, pp. 293--311, 1988.

\bibitem{Kundur1994}
P.~Kundur, \emph{{Power System Stability and Control}}, 1st~ed., N.~J. Bau and
  M.~G. Lauby, Eds.\hskip 1em plus 0.5em minus 0.4em\relax McGraw-Hill
  Professional, 1994.

\bibitem{Borsche2016Diss}
T.~S. Borsche, ``Impact of demand and storage control on power system operation
  and dynamics,'' Ph.D. dissertation, ETH Z\"urich, Feb 2016.

\end{thebibliography}

\end{document}